\documentclass[12pt]{article}
\usepackage{amsmath,amssymb}
\usepackage{authblk}
\usepackage{color}
\usepackage{mathrsfs}
\usepackage{mathptmx}
\usepackage{verbatim}
\usepackage{epsfig}
\usepackage{CJK}
\usepackage{bm}
\usepackage{amsfonts}
\usepackage{amsthm}
\input{epsf.sty}
\usepackage{epsfig}
\usepackage{setspace}

\def\<{\langle}
\def\>{\rangle}

\def\be{\begin{equation}}
\def\ee{\end{equation}}
\def\ba{\begin{array}}
\def\ea{\end{array}}

\newtheorem{theorem}{Theorem}[section]
\newtheorem*{thmA}{Theorem A}

\newtheorem*{thmC}{Theorem C}
\newtheorem{lemma}{Lemma}[section]
\newtheorem{proposition}{Proposition}[section]

\theoremstyle{definition}
\newtheorem{remark}{Remark}[section]
\newtheorem{example}{Example}[section]
\newtheorem{definition}{Definition}[section]
\newtheorem{problem}{Problem}
\numberwithin{equation}{section}

\usepackage{amsmath}
\usepackage{amsthm}
\usepackage{amsfonts}
\usepackage{pifont}
\usepackage{color}
\usepackage{bbm}
\usepackage{CJK}
\usepackage{mathrsfs}
\usepackage{fancyhdr}
\usepackage{graphicx}
\usepackage{pdfsync}
\usepackage{psfrag}
\usepackage{time}
\usepackage{txfonts}

\usepackage{stmaryrd}
\usepackage{mathrsfs}
\usepackage{txfonts}
\usepackage{amsfonts}

\usepackage{indentfirst}

\def\be{\begin{equation}}
\def\ee{\end{equation}}
\def\br{\begin{eqnarray}}
\def\er{\end{eqnarray}}

\title{On the Ergodic theory of the           Generalized incompressible flow
\thanks
{Supported by   NSFC11271076 and  NSFC11121101   }}

\author{$\mbox{Cheng Yang\thanks{
E-mail:11110180006@fudan.edu.cn} \;  and \; Xiaoping Yuan\thanks{
E-mail:xpyuan@fudan.edu.cn   }}$
\\
School of Mathematical Sciences,\ \  Fudan University\\
Shanghai  200433,\ \   P.R.China
\\
}

\begin{document}

\maketitle

\date{}

\begin{quote}
\small {\bf Abstract.}  In the  study of the variational problem related to the incompressible fluid, Brenier brought the concept of generalized flows and showed that the generalized incompressible flow (GIF) is deeply related to the study of the incompressible Euler equations. It turns out that the classical ergodic theory does not apply to  GIF which, regarded as a very special stochastic process, is neither stationary nor Markov. In the present paper we develop a generalized ergodic theory for GIF.

\end{quote}

\tableofcontents

\section{Introduction}
In his study of the variational problem of ideal fluid motions, Brenier brought the celebrated idea of generalized incompressible flows (GIF). A GIF is a generalization of the classical incompressible flow and can be regarded as a very special stochastic process which is neither stationary nor Markov. So the classical ergodic theory can not apply to GIFs. Naturally, some kind of generalized ergodic theory should be developed for GIF.


\subsection{Backgrounds on the incompressible fluids and generalized flows}
\par
In his paper \cite{arn} published in 1966, V. I. Arnold studied the geometric approach to the incompressible Euler equations. This work brings many new ideas into the research of the ideal incompressible fluid. See \cite{arn2}, \cite{khe} for the details.
\par
Let us recall that the incompressible Euler equations on  a compact, oriented Riemannian manifold $M$ without boundary. Denote the Riemannian metric by $(\;,\;)$ and the volume form by $\mu$. The Euler equations of incompressible flows of homogenous fluid on the manifold M are
\begin{equation}\label{e11}
 \left\{
  \begin{array}{ll}
    div_{\mu} v= 0 , & \hbox{} \\
     \frac{Dv}{Dt}= -\nabla p, & \hbox{} \\
    v\mid_{t=0}=v_0,  & \hbox{}
  \end{array}
\right.
\end{equation}
where $\frac{Dv}{Dt}=\frac{\partial v}{\partial t}+\nabla_{v}v$, $\nabla_{v}v$ is the covariant derivative of the velocity field $v(x,t)$ along itself, and $p(x,t)$ is the scalar pressure field which can be determined by the condition $div_{\mu}(\frac{\partial v}{\partial t})=0$ up to a time-dependent constant. The first equation shows that $v$ is divergence-free vector field, i.e. the flow of it preserve the volume form $\mu$.

\par
One can study (\ref{e11})  via those theories of partial differential equations. Refer to \cite{ber}, \cite{lio} for example. An alternate way in the  study
of  (\ref{e11}) is Lagrangian formulation. When one's attention is focused on the trajectories of particles,  the incompressible fluid can be considered as a map $t\mapsto g(t,\cdot)$, here $g(t,\cdot)$ is a volume preserving diffeomorphism. So we can say that the configuration space is $\text{SDiff}(M)$, here $\text{SDiff}(M)$ denotes the group of volume-preserving diffeomorphisms from $M$ to itself.
The group $\text{SDiff}(M)\subseteq L^{2}(M;\mathbb{R}^{n})$ can be seen as an infinite-dimensional Hilbert manifold and the energy
\begin{equation}\label{e14}
E=\frac{1}{2} \int_M(v,v)\;dvol_M
\end{equation}
introduces a right-invariant weak-Riemannian metric on $\text{SDiff}(M)$. Actually, this Riemannian structure is inherited from the Hilbert space $L^{2}(M;\mathbb{R}^{n})$. Arnold \cite{arn} pointed out that the incompressible fluid is a geodesic on the group $\text{SDiff}(M)$. Therefore the study of Euler equation can be converted into the following variational problem:\\
\begin{problem}   \label{p11}
Given $g, h\in \text{SDiff}(M)$, find a smooth curve $\{g_{t}\}_{t\geq0}\subset \text{SDiff}(M)$ minimizes the energy
\begin{equation}\label{e14}
  E(g)=\int^{1}_{0}\frac{1}{2}\|\dot{g}_{t}\|^{2}_{L^{2}(M;\mathbb{R}^{n})}\;dt=\int^{1}_{0}\int_{M}\frac{1}{2}(\frac{\partial g(t,x)}{\partial t},\frac{\partial g(t,x)}{\partial t})\;dtdx,
\end{equation}
among all curves satisfying $g_{0}=g,g_{1}=h$.
\end{problem}
\par
An local existence and uniqueness theorem for this problem has been proved by Ebin and Marsden \cite{ebi} if $g$ and $h$ are sufficiently close in a sufficiently high order Sobolev norm.
\par
In 1985, however, Shnirelman \cite{shn1} introduced a diffeomorphism $g\in \text{SDiff}(K^{3})$ such that there is no minimizing curve connecting the identity map to $g$, here $K^3$ is the cube $[0,1]^3$. Shnirelman's counterexample was constructed via a special class of diffeomorphisms $g(x,z)=(h(x),z)$, here $g\in \text{SDiff}(K^{3}),\;h\in \text{SDiff}(K^{2})$. Because the fluid motions in two dimensional plane have more constraints, a minimizing path connecting the identity map to $g$ would have to use the third dimension. Then we compress the fluid movements in the third dimension, which leads to the new movement having less kinetic energy than before.
\par
In 1989, Y.Brenier \cite{bre} considered the concept of generalized flows, this concept is closely related to the  Young's measure introduced by L. C. Young. Then he studied the generalized variational problem, the relation between this variational problem and problem \ref{p11} is like the Kantorovitch problem and Monge problem in optimal transport. Next we introduce Brenier's work briefly.
\par
Because of the nonlinear constraint, the variational problem 1 is very hard to answer. In the spirit of Young's work in calculus of variations, Brenier brought the concept of generalized flows which is a probability measure on the product space $\Omega=M^{[0,T]}$,  $\Omega$ can be viewed as the space of all the paths $t\in[0,T]\rightarrow z(t)\in M$ on M. Then he proved the existence of the minimizing solution of the following variational problem \cite{bre}:
\begin{problem}   \label{p13}
Given $h\in MP(M)=\{\gamma:M\rightarrow M\mid$ if $Y$ is a measurable subset of $M$, then $\gamma^{-1}(Y)$ is measurable and $\text{meas}\;(\gamma^{-1}(Y))=\text{meas}\; (Y)\}$, find a generalized incompressible flow $q$ minimizing the energy
\begin{equation}\label{e19}
  E(q)=\int\rho(z(0))\int\frac{1}{2}\parallel z'(t)\parallel^2\;dt\;q(dz),
\end{equation}
among all the generalized incompressible flows which reach the final configuration $\eta_h$, i.e. satisfies
\begin{equation}\label{e20}
  \int_{\Omega}f(z(0),z(t))\;q(dz)=\int_{M\times M}f(x,y)\;\eta_h(dx,dy),\;\forall f\in C(M\times M),
\end{equation}
where a final configuration $\eta_h(dx,dy)$ is a doubly stochastic probability measure $\eta_h$ defined by
\begin{equation}\label{e18}
  \eta_h(dx,dy)=\delta(y-h(x))dx.
\end{equation}
\end{problem}
Then Brenier proved that under a natural restriction on the time scale, any classical solution to the Euler equations satisfies the generalized variational problem and is the unique minimizing solution. So the new concept of generalized flows did not miss the classical solution. Brenier also gives a formal relation between his generalized flows and Diperna and Majda's measure valued solutions \cite{dip}.
\par
The idea of generalized flows means that each fluid particle may split and move in different direction with certain probability, even cross each other. Using the notion of generalized flows, Brenier proved the existence of solutions for the relaxed short path problem \cite{bre2}. He also showed the notion of generalized flow is related to the hydrostatic approximation of the Euler equation \cite{bre3}.
\par
The technique of GIFs is also proved to be a powerful tool for studying the geometry of the volume-preserving diffeomorphism group $\text{SDiff}(M)$. In \cite{shn2}, Shnirelman obtained the accurate estimate of the diameter and the degree of regularity of the volume-preserving diffeomorphism group $\text{SDiff}(M)$ by using the idea of GIFs.
\par
 Now we can see that the notion of generalized incompressible flows (GIFs) not only is deeply related to incompressible Euler equation, but also provide a useful tool to investigate the geometry of the volume-preserving diffeomorphism group $\text{SDiff}(M)$. Its dynamic property may shed some lights on the study of incompressible fluid. A natural question is how to develop a generalized ergodic theory for GIFs? More specific, we may ask following questions: what are the recurrence and ergodic properties of GIFs? Are there any ergodic theorems of GIFs? What is the structure of the set of all GIFs? We will answer these questions in our paper. As we will mention in section 2, the set of GIFs is indeed larger than the set of invariant measures on $\Omega=M^{[-\infty,+\infty]}$. So there are many differences between the ergodic theory of GIFs and the classical one.
\subsection{Organization and main results.}
\par
The present paper is organized as follows:
\par
In \S 2, after defining the concept of GIFs, we will study the recurrence property of GIFs. We will show that GIFs only have the recurrence property in a weak sense. See Theorem \ref{t21}.
We also construct a counterexample to show that GIFs don't satisfy the strong version of the recurrence property.
\par
In \S 3, we will define the concept of the ergodicity for GIFs. We will prove that ergodic GIFs include the classical ergodic flows by Proposition 3.1.
Also, we will give some equivalent definitions of the ergodicity in Proposition 3.3.
\par
In \S 4, we will prove several ergodic theorems for GIFs. Our method is non-traditional, and it is closer to Bourgain's approach to the ergodic theorem \cite{bou4}. Bourgain developed his method in the study of the extension of Birkhoff's pointwise ergodic theorem to certain subset of integrals and its application to analytic number theory \cite{bou1}, \cite{bou2}, \cite{bou3}. First we will prove a maximal ergodic theorem (Theorem \ref{t41}) by the harmonic analysis method, then we will prove the $L^1$ ergodic theorem (Theorem \ref{t43}). Finally combining these two theorems we will prove a pointwise ergodic theorem (Theorem \ref{t46}) for GIFs despite the fact that the GIFs only have weak recurrence property.

\par
In \S 5, we will study the set $\text{GIF}\,(\Omega)$ of all GIFs. Like the set of invariant measures on a compact metric space, the $\text{GIF}\,(\Omega)$ is a convex and compact set. The difference is that the extreme points of the set are not ergodic GIFs. We will prove that the extreme points of the set are weak ergodic GIFs in Theorem \ref{t54}.

\bigskip

{\bf  Acknowledgements.}
The first author is grateful to Boris Khesin for comments and fruitful discussions.

\bigskip

\section{Generalized incompressible flows and recurrence property  }
First we give the explicit definitions of the generalized incompressible flows. Let X be a compact metric space, $\Omega=X^\mathbb{R}$ is the path space on X.
\begin{definition} \textbf{(Generalized incompressible flow)}  \label{d21}
 A probability measure q on $\Omega$ is called a generalized incompressible flow (GIF) if it satisfies the following incompressible condition:
 \begin{equation}\label{e21}
  \int_{\Omega}f(z(t_{0}))\;q(dz)=\int_{\Omega}f(z(0))\;q(dz), \forall f(x)\in C(X), \forall t_{0}\in \mathbb{R},
\end{equation}
\end{definition}
\begin{remark}   \label{r21}
 (a) By the definition, since continuous functions on $\Omega$ of type $f(z(t))$ are a subset of all continuous functions on $\Omega$, we cannot claim that q is an invariant measure of the flow $\Gamma^s$ on $\Omega$:
  \begin{eqnarray*}\label{e22}
  \Gamma^s:&&\Omega\longrightarrow\Omega \\
   &&z(t)\longmapsto z(t+s).
\end{eqnarray*}
 However, for an invariant measure of $\Gamma^s$, it must be a GIF. The set consisting of GIFs on $\Omega$ 
 is larger than the set of invariant measures on $\Omega$.
\par
(b) Let $E\subset X$ be a Borel set of X. We have $\int_{\Omega}\chi_{E}(z(t_{0}))\;q(dz)=\int_{\Omega}\chi_{E}(z(0))\;q(dz)$, hence we can define $\mu(E)=\int_{\Omega}\chi_{E}(z(0))\;q(dz)$. This gives a probability measure on $X$, such that
 \begin{equation}\label{e23}
  \int_{\Omega}f(z(t_{0}))\;q(dz)=\int_{\Omega}f(z(0))\;q(dz)=\int_{X}f(x)\;\mu(dx), \forall f(x)\in C(X), \forall t_{0}\in \mathbb{R}.
\end{equation}
\par
(c) As mentioned by Shnirelman in \cite{shn2}, GIFs can be regarded as a stochastic process:
\begin{definition}\label{d210}
Let $(\Omega, \mathcal{B}, \mathbb{P})$ be a probability space which can be seen as the space of fluid particles. The GIF on a manifold $M$ is a stochastic process:
\begin{eqnarray*}\label{e210}
  x:&&\Omega\times\mathbb{R}\longrightarrow M \\
   &&(\omega,\;t)\longmapsto x(\omega,\;t)
\end{eqnarray*}
satisfying the incompressible condition:
$$
\mathbb{P}\{\omega\mid x(\omega,\;t_0)\in X\}=\text{Vol}_M(X)
$$
for all $t_0\in\mathbb{R}$ and $X\subset M$, where $\text{Vol}_M$ is the volume form on $M$.
\end{definition}
We can see that a GIF is a very special stochastic process, it is neither stationary nor Markov.
 \end{remark}
Let $(\Omega,q)$ be a GIF, it has the recurrence property in the following sense:
\begin{theorem}   \label{t21}
 Let $(\Omega,q)$ be a GIF, and assume we have a Borel set $E\in \mathcal{B}_{X}$ with a positive measure $\mu(E)>0$, here $\mu$ is a probability measure on X satisfies (\ref{e23}). Then for the subset  $\Phi=\{z(t)\in \Omega\mid z(0)\in E\}$ of all the paths start from $E$ , there exists $s\in \mathbb{R}^{+}$ such that $q(\Phi\cap\Gamma^{-s}\Phi)>0$.
\end{theorem}
\noindent\textbf{Proof.} Let us consider the set $\{\Gamma^{-s}\Phi\mid s>0\}$. If we have
 \begin{equation}\label{e231}
  q(\Gamma^{-s}\Phi\cap\Gamma^{-t}\Phi)=0,\forall s\neq t>0,
\end{equation}
then $q(\cup_{s>0}\Gamma^{-s}\Phi)=\infty$. This is because $\forall s>0$,
\begin{eqnarray*}
 q(\Gamma^{-s}\Phi) &=& \int_{\Omega}\chi_{\Gamma^{-s}\Phi}(z(t))\;q(dz)=\int_{\Omega}\chi_{\Phi}(z(t+s))\;q(dz)\\
 &=& \int_{\Omega}\chi_{E}(z(s))\;q(dz)=\int_{\Omega}\chi_{E}(z(0))\;q(dz)=\mu(E)>0.
\end{eqnarray*}
Therefore there exist $s<t$ such that $q(\Gamma^{-s}\Phi\cap\Gamma^{-t}\Phi)>0$, i.e. $q(\Phi\cap\Gamma^{-(t-s)}\Phi)>0$.
\begin{flushright}
  $\Box$
\end{flushright}
We may guess the strong version of the recurrence property also holds:\\
\noindent\textbf{Conjecture.} Let $(\Omega,q)$ be a GIF, and we have a Borel set $E\in \mathcal{B}_{X}$ with a positive measure $\mu(E)>0$, here $\mu$ is a probability measure on X satisfying (\ref{e23}). Then for $\Phi=\{z(t)\in \Omega\mid z(0)\in E\}$, there exists $\Theta\subset\Phi$ such that $q(\Theta)=q(\Phi)$ and $\Theta=\{z(t)\in\Omega\mid\exists t_{1}<\cdots<t_{n}<\cdots$, such that $z(t_{i})\in E\}$.
\par
Unfortunately, this conjecture is not true. Here we give a counterexample which is a GIF, but for $q$ almost every $z(t)\in\Omega$, the flow $z(t)$ does not have the recurrence property.
\begin{example}   \label{ex21}
 Let X be the unit circle $S^{1}$. Considering the following flow $z(t)$ on $S^{1}$:
  \begin{equation}\label{e24}
  z_x(t)=\left\{
     \begin{array}{ll}
       x+t(mod\;2\pi) & ,t\in[0,\frac{\pi}{2}] \\
       x+\frac{\pi}{2}(mod\;2\pi) & ,t\in(\frac{\pi}{2},+\infty)
     \end{array}
   \right.
\end{equation}
We will study the generalized flow $q$ which equidistributes on the set of the flows defined above. It is easy to check that
\begin{eqnarray*}
  \int_{\Omega}\chi_{(s_{1},s_{2}]}(z(t))\;q(dz)=\int_{\Omega}\chi_{(s_{1},s_{2}]}(z(0))\;q(dz)
  =\int_{X}\chi_{(s_{1},s_{2}]}(x)\;dx=\frac{s_{2}-s_{1}}{2\pi},
\end{eqnarray*}
hence $q$ is a GIF. Since the flow defined above stops at $t=\frac{\pi}{2}$, so for $q$ almost every $z(t)$, the flow does not have the recurrence property.

\end{example}
\begin{remark}\label{r22}
From this example we can see that the GIF only has the recurrence property in a very weak sense. Even so, in the section 4, we will prove a pointwise ergodic theorem of the GIF.
\end{remark}

\section{Ergodicity}
In this section, we will discuss the ergodicity of the GIF. For a good introduction to classical ergodic theory, one can refer to \cite{ein}, \cite{wal}. First, we give the definition of ergodicity.

\begin{definition}   \textbf{(Ergodic GIF)}\label{d31}
A GIF $q$ is called ergodic if for any Borel set $B\in \mathcal{B}_{\Omega}$ which satisfies $\Gamma^{-t}(B)=B, \forall t\in \mathbb{R}$, we have $q(B)=0$ or $1$.
\end{definition}
\par
Next we show that the GIFs induced by classical ergodic flows are ergodic GIFs.
\par
Let $j(t)$ be a classical ergodic flow on $(X,\mu)$. For $\mu-$almost all $x\in X$, we can find a unique $t_{0}$ such that $x=j(t_{0})$, then we can define $G(t,x)=j(t_{0}+t)$ for $\mu-$almost all $x\in X$. So a generalized flow $q$ induced by $j(t)$ is defined as follows:
\begin{equation}\label{e31}
  \int_{\Omega}F(z(t))\;q(dz)=\int_{X}F(x\mapsto G(t,x))\;\mu(dx), \forall F(z)\in C(\Omega).
\end{equation}
We claim that $q$ is a GIF. Indeed, $\forall f(x)\in C(X)$, we have
\begin{equation}\label{e32}
  \int_{\Omega}f(z(t_{0}))\;q(dz)=\int_{X}f(G(t_{0},x))\;dx=\int_{X}f(x)\;\mu(dx),\forall t_0\in \mathbb{R},
\end{equation}
the second equality follows from that $\mu$ is the invariant measure for the flow $j(t)$. Now we prove that q is an ergodic GIF.
\begin{proposition}   \label{p31}
 The GIF $q$ defined by (\ref{e31}) above is ergodic.
\end{proposition}
\noindent\textbf{Proof.} Let $Q=\{u(t)=j(t_{0}+t),\forall t_{0}\in \mathbb{R}\}$. Then, by the definition of q (see (\ref{e31})), we get $q(Q)=1$. Let $B$ be a Borel set of $\Omega$ with $\Gamma^{-t}(B)=B, \forall t\in \mathbb{R}$.
\par
If $B\cap Q=\emptyset$, then $q(B)=0$. If $B\cap Q\neq\emptyset$, then we can find $j(t_{0}+t)\in B\cap Q$. Since $\Gamma^{-t_{0}}(j(t_{0}+t))=j(t)$ and $B$ is invariant, we get $\{\Gamma^{-t_{0}}(j(t))\}\subset Q\cap B$, therefore $q(B)=1$. We see that $q(B)$ is either 0 or 1, so $q$ is an ergodic GIF.
\begin{flushright}
  $\Box$
\end{flushright}
\par
Before continuing our discussion, we first prove the following Lemma. We will always assume that the map $(z(t),s)\mapsto \Gamma^{s}z(t)$ is measurable.
\begin{lemma}   \label{p32}
 Let $(\Omega,q)$ be a GIF. Then for $B\in \mathcal{B}_{\Omega}$, the following statements are equivalent:\\
 $(1)\;q(\Gamma^{s}B\triangle B)=0, \forall s\in \mathbb{R}$.\\
 $(2)$ There exists $F\in \mathcal{B}_{\Omega}$ such that $q(F\triangle B)=0$ and $\Gamma^{s}F=F,\forall s\in \mathbb{R}$.
\end{lemma}
\noindent\textbf{Proof.} $(2)\Rightarrow (1)$ is obvious. Now we prove $(1)\Rightarrow (2)$. First taking the rational number set $\mathbb{Q}$, define $F=\cap_{s\in \mathbb{Q}}\Gamma^{s}B\subseteq B$. Then $\Gamma^{s}F=F,\forall s\in \mathbb{Q}$ and $q(B\triangle F)=q(B\backslash F)=q(\cup_{s\in \mathbb{Q}}B\backslash\Gamma^{s}B)=0$.
\par
Suppose $B\in \mathcal{B}_{\Omega}$ satisfies $(1)$. Then by the discussion above we can assume that $\Gamma^{s}B=B,\forall s\in \mathbb{Q}$. Define $B_{z}=\{s\in \mathbb{R}\mid\Gamma^{s}z \in B\}\subseteq \mathbb{R}$. By the assumption that $(z(t),s)\mapsto \Gamma^{s}z(t)$ is measurable, we know $B_{z}$ is measurable for each $z\in\Omega$.
\par
Because $\Gamma^{s}B=B,\forall s\in \mathbb{Q}$, we get $s+B_{z}=B_{z},\forall s\in \mathbb{Q}$. Using Proposition $8.6$ in \cite{ein}, we know that either $m(B_{z})=0$ or $m(\mathbb{R}\backslash B_{z})=0$, here m is the Lebesgue measure on $\mathbb{R}$. Hence we can define $F=\{z\in\Omega\mid m(B_{z})>0\}=\{z\in\Omega\mid m(\mathbb{R}\backslash B_{z})=0\}$. If $z\in F$ and $s\in \mathbb{R}$, then $B_{\Gamma^{s}z}=\{r\in\mathbb{R}\mid\Gamma^{r}\Gamma^{s}z\in B\}=B_{z}-s$. Using Proposition $8.6$ in \cite{ein} again, we see that $m(B_{z})>0$ implies $m(B_{\Gamma^{s}})>0$, i.e. $\Gamma^{s}z\in F$. So $\Gamma^{s}F=F,\forall s\in \mathbb{R}$.
\par
Last, we prove that $F$ is measurable and $q(B\triangle F)=0$. Let $I\subseteq \mathbb{R}$ be a measurable set of positive measure and define
\begin{equation}\label{e32}
  f(z)=\frac{1}{m(I)}\int_{I}\chi_{B}(\Gamma^{s}z)\;dm(s).
\end{equation}
We see that $f(z)=\chi_{F}$. Hence by Fubini's theorem, $F$ is measurable. Noticing that $\chi_{B}(\Gamma^{s}z)=\chi_{\Gamma^{-s}B(z)}$, we get
\begin{eqnarray*}
 q(B\triangle F) &=& \int_{\Omega}(\chi_{B}+\chi_{F}-2\chi_{B}\chi{F})\;q(dz)\\
 &=&\frac{1}{m(I)}\int_{I}\int_{\Omega}(\chi_{B}+\chi_{\Gamma^{-s}B}-2\chi_{B}\chi_{\Gamma^{-s}B})q(dz)\;m(ds)\\
 &=&\frac{1}{m(I)}\int_{I}q(B\triangle\Gamma^{-s}B)\;m(ds)=0.
\end{eqnarray*}
\begin{flushright}
  $\Box$
\end{flushright}
\par
Now we could give several equivalent way to define the ergodicity.
\begin{proposition}   \label{p33}
 Let $(\Omega,q)$ be a GIF, then the following statements are equivalent:\\
 $(1)$ q is ergodic.\\
 $(2)$ If $\Phi\in\mathcal{B}_{\Omega}$ with $q(\Gamma^{s}\Phi\triangle \Phi)=0,\forall s\in\mathbb{R}$, we have $q(\Phi)=0$ or $1$.\\
 $(3)$ If $f$ is $q-$measurable and $f(\Gamma^{s}z)=f(z),\forall s\in\mathbb{R},\forall z\in\Omega$, then for $q-$almost every $z$, $f$ is constant.
\end{proposition}
\noindent\textbf{Proof.}$(1)\Rightarrow (2).$ If $\Phi\in\mathcal{B}_{\Omega}$ satisfies $q(\Gamma^{s}\Phi\triangle \Phi)=0,\forall s\in\mathbb{R}$, then by Lemma \ref{p32}, we get $\Psi\in \mathcal{B}_{\Omega}$ such that $q(\Phi\triangle\Psi)=0$ and $\Gamma^{s}\Psi=\Psi,\forall s\in \mathbb{R}$. By definition \ref{d31}, we know that $q(\Psi)=0$ or 1, i.e. $q(\Phi)=0$ or 1.\\
$(2)\Rightarrow (1).$ If $\Gamma^{s}\Phi=\Phi,\forall s\in \mathbb{R}$, so $q(\Gamma^{s}\Phi\triangle \Phi)=0,\forall s\in\mathbb{R}$. By $(2)$, we have $q(\Phi)=0$ or $1$ which proves the ergodicity.\\
$(1)\Rightarrow (3).$ We can assume that $f$ is real-valued function. If $q$ is ergodic, we define
\begin{equation}\label{e33}
  Z(k,n)=\{z\in \Omega\mid\frac{k}{2^n}\leq f(z)<\frac{k+1}{2^n}\}=f^{-1}([\frac{k}{2^n},\frac{k+1}{2^n}]),
\end{equation}
here $k\in \mathbb{Z}$ and $n\in\mathbb{N}^{+}$. We have $\Gamma^{-s}Z(k,n)\triangle Z(k,n)\subseteq \{z\mid f(\Gamma^{s}z)\neq f(z)\}$, so $q(\Gamma^{-s}Z(k,n)\triangle Z(k,n))=0, \forall s\in\mathbb{R}$ which means $q(Z(k,n))=0$ or 1 by the ergodicity of q.
\par
For fixed n,$\cup_{k\in\mathbb{Z}}Z(k,n)=\Omega$ is a disjoint union and so there exists a unique $k_{n}$ with $q(Z(k_{n},n))=1$. Let $Q=\cap_{n=1}^{\infty}Z(k_{n},n)$, then $q(Q)=1$ and $f$ is constant on Q. \\
$(3)\Rightarrow (1).$ Suppose $\Gamma^{-s}\Phi=\Phi,\forall s\in\mathbb{R}$, then $\chi_{\Phi}$ is measurable and $\chi_{\Phi}(\Gamma^{s}z)=\chi_{\Phi}(z), \forall z\in\Omega$. So $\chi_{\Phi}$ is constant $q-$almost every $z$, i.e. $q(\Phi)=0$ or 1.

\begin{flushright}
  $\Box$
\end{flushright}

\section{The ergodic theorems of GIFs }
Let X be a compact metric space, path space $\Omega=X^{[0,+\infty)}$. In our following discussion, we will assume that the map
\begin{eqnarray*}\
  \Gamma:(\Omega \times \mathbb{R},q \times m)&\longrightarrow &(\Omega,q) \\
   (z(t),t_{0})&\longmapsto &z(t+t_{0})
\end{eqnarray*}
is measurable. First we prove the maximal ergodic theorem.
\begin{theorem}   \textbf{(Maximal ergodic theorem)}\label{t41}
 Let $(\Omega,q)$ be a GIF, for $f(x)\in L^1(X)$ and $\alpha >0$, we define the maximal function
 \begin{equation}\label{e41}
  f^{*}(z)=\sup_{T>0}\frac{1}{T}\int_{0}^{T}f(z(t))\;dt,
\end{equation}
and the set
 \begin{equation}\label{e42}
  E_{\alpha}^{f}=\{z\in\Omega\mid f^{*}(z)>\alpha\},
\end{equation}
then we have
 \begin{equation}\label{e43}
  \alpha q(E^{f}_{\alpha})\leq\parallel f\parallel_{L^1(X)}.
\end{equation}

\end{theorem}
\par
Before we prove this theorem, we introduce the finite Vitali covering lemma on $\mathbb{R}$,
\begin{lemma}  \textbf{(Finite Vitali covering lemma)} \label{t42}
 For any collection of intervals $I_{1}=[a_1,a_1+l_1],\ldots,I_K=[a_K,a_K+l_K]$ in $\mathbb{R}$, there is a disjoint subcollection $I_{j_1},\ldots,I_{j_k}$ such that
 \begin{equation}\label{e44}
  I_1\cup\cdots\cup I_K\subseteq\bigcup_{m=1}^{k}[a_{j_m}-l_{j_m},a_{j_m}+2l_{j_m}].
\end{equation}

\end{lemma}
\par
Use the finite Vitali covering lemma, we prove the following lemma.
\begin{lemma} \label{t43}
$ \forall g\in L^1(\mathbb{R})$ and $\alpha>0$, define the maximal function
\begin{equation}\label{e45}
  g^{*}(a)=\sup_{T>0}\frac{1}{T}\int_{0}^{T}g(a+t)\;dt,
\end{equation}
and the set
\begin{equation}\label{e46}
  E_{\alpha}^{g}=\{a\in\mathbb{R}\mid g^{*}(a)>\alpha\}.
\end{equation}
Then we have
\begin{equation}\label{e47}
  \alpha \mid E^{g}_{\alpha}\mid\leq 3\parallel g\parallel_{L^1(\mathbb{R})}=3\int_{-\infty}^{+\infty}g(t)\;dt.
\end{equation}

\end{lemma}
\noindent\textbf{Proof.} For any $a\in E_{\alpha}^{g}$, choosing $l(a)$ such that $\frac{1}{l(a)}\int_{0}^{l(a)}g(a+t)\;dt>\alpha$. Since $\mathbb{R}$ is second countable, we can write the set $O=\cup_{a\in E_{\alpha}^{g}}[a,a+l(a)]$ as a union of countable many sets $O=\bigcup_{i=1}^{\infty}[a_i,a_i+l(a_i)]$.\\
Fix some $K\geqslant1$, and using the covering lemma, we get the subcollection $I_{i_1},\ldots,I_{i_k}$. Since they are disjoint, it follows that
\begin{equation}\label{e48}
  \sum_{j=1}^{k}\int_{I_{i_j}}g(t)\;dt\leqslant \parallel g\parallel_{L^1(\mathbb{R})},
\end{equation}
where the left hand side equals to
\begin{equation}\label{e49}
  \sum_{j=1}^{k}l(a_{i_j})\frac{1}{l(a_{i_j})}\int_{0}^{l(a_{i_j})}g(a_{i_j}+t)\;dt>\sum_{j=1}^{k}l(a_{i_j})\alpha,
\end{equation}
by the choice of $l(a)$. However, since
\begin{equation}\label{e50}
 I_1\cup\cdots\cup I_K\subseteq\bigcup_{j=1}^{k}[a_{i_j}-l(a_{i_j}),a_{i_j}+2l(a_{i_j})],
\end{equation}
we have
\begin{equation}\label{e51}
 \mid I_1\cup\cdots\cup I_K\mid\leq 3\sum_{j=1}^{k}l(a_{i_j}).
\end{equation}
By (\ref{e49}) and (\ref{e51}), we have
\begin{equation}\label{e52}
 \alpha\mid I_1\cup\cdots\cup I_K\mid\leq 3\parallel g\parallel_{L^1(\mathbb{R})}, \text{ for any } K\geq1,
\end{equation}
which gives
\begin{equation}\label{e52}
 \alpha m(O)\leq 3\parallel g\parallel_{L^1(\mathbb{R})}.
\end{equation}
Then by $E_{\alpha}^{g}\subseteq O$, we get the lemma.
\begin{flushright}
  $\Box$
\end{flushright}

\noindent\textbf{Proof of Theorem \ref{t41}.} We define
\begin{equation}\label{e53}
 f^{*}_{M}(z)=\sup_{0\leq T\leq M}\frac{1}{T}\int_{0}^{T}f(z(t))\;dt,
\end{equation}
and the set
\begin{equation}\label{e54}
 E_{\alpha,M}^{f}=\{z\in\Omega\mid f^{*}_{M}(z)>\alpha\}.
\end{equation}
Given $J>0$, let
\begin{equation}\label{e55}
  g(t)=\left\{
     \begin{array}{ll}
       f(z(t)) & ,t\in[0,J] \\
       0 & ,otherwise.
     \end{array}
   \right.
\end{equation}
Similarly, we can define
\begin{equation}\label{e56}
 g^{*}_{M}(a)=\sup_{0\leq T\leq M}\frac{1}{T}\int_{0}^{T}g(a+t))\;dt.
\end{equation}
If $0\leq a\leq J-M$ and $0\leq t\leq M$, we have
\begin{equation}\label{e57}
 g(a+t)=f(z(a+t)),
\end{equation}
so
\begin{equation}\label{e58}
 g^{*}_{M}(a)=f^{*}_{M}(z(a+t)), \ \ 0\leq a\leq J-M.
\end{equation}
According to Lemma \ref{t43}, $\forall z\in\Omega,\ \alpha >0$, we have
\begin{eqnarray*}\
  \alpha\mid\{a\in[0,J-M]\mid g^{*}_{M}(a)>\alpha\}\mid &\leq& \alpha\mid\{a\in[0,J]\mid g^{*}_{M}(a)>\alpha\}\mid\\
  &\leq& \alpha\mid\{a\in\mathbb{R}\mid g^{*}(a)>\alpha\}\mid\\
  &\leq& 3\parallel g \parallel_{L^1(\mathbb{R})},
\end{eqnarray*}
or
\begin{eqnarray*}\
  \alpha\int_{0}^{J-M}\chi_{E^{f}_{\alpha,M}}(z(a+t))\;da &=& \alpha\mid\{a\in[0,J-M]\mid f^{*}_{M}(z(a+t))>\alpha\}\mid\\
  &\leq& 3\int_{0}^{J}f(z(t))\;dt.
\end{eqnarray*}
Integrating both sides of the inequality over $\Omega$ and using the Fubini's Theorem, we can get
\begin{equation}\label{e59}
 (J-M)\alpha q(E^{f}_{\alpha,M})\leq 3J\parallel f\parallel_{L^1_{q}(\Omega)}.
\end{equation}
Letting $J\rightarrow\infty$ gives
\begin{equation}\label{e60}
 \alpha q(E^{f}_{\alpha,M})\leq 3\parallel f\parallel_{L^1_{q}(\Omega)}=3\parallel f\parallel_{L^1_{\mu}(X)}.
\end{equation}
Finally, letting $M\rightarrow\infty$, we get theorem \ref{t41}.
\begin{flushright}
  $\Box$
\end{flushright}
\par
Recall that we have defined the flow $\Gamma^s$ on $\Omega$ in section 2:
\begin{eqnarray*}\label{e22}
  \Gamma^s:&&\Omega\longrightarrow\Omega \\
   &&z(t)\longmapsto z(t+s).
\end{eqnarray*}
Then we can define an operator $U_{s}$ as follows
\begin{eqnarray*}\label{e22}
  U_{s}:&&L^2_q(\Omega)\longrightarrow L^2_q(\Omega) \\
   &&f(z(t))\longmapsto f(\Gamma^s(z(t))).
\end{eqnarray*}
Let $f(x)\in L^2_{\mu}(X)$, then
\begin{equation}\label{e61}
\int_{\Omega}\mid f(z(t_{0}))\mid^2\;q(dz)=\int_{\Omega}\mid f(z(t_{0}+s))\mid^2\;q(dz).
\end{equation}
Actually, if $f(x)\in L^2_{\mu}(X)$, then $f(z(t_{0}))\in L^2_q(\Omega)$. Hence $\{f(z(t_{0}))\mid f(x)\in L^2_{\mu}(X),t_{0}\in\mathbb{R}\}$ can generate a linear subspace of Hilbert space $L^2_q(\Omega)$. Denote by $H$ the $L^2_{q}-$ closure of the subspace. So $H$ is a Hilbert space. Because
\begin{eqnarray*}
\int_{\Omega}\mid af(z(t_{0}))+bg(z(t_1))\mid^2 \;q(dz) &=& \int_X\mid af(x)+bg(x))\mid^2 \;\mu(dx)\\
 &=& \int_{\Omega}\mid af(z(t_{0}+s))+bg(z_1+s)\mid^2 \;q(dz).
\end{eqnarray*}
and
\begin{eqnarray*}
\int_{\Omega}\mid f(z(t_{0}))\mid^2 \;q(dz) &=& \lim_{n\rightarrow\infty}\int_{\Omega}\mid f_n(z(t_0))\mid^2 \;q(dz)\\
 &=& \lim_{n\rightarrow\infty}\int_{\Omega}\mid f_n(z(t_0+s))\mid^2 \;q(dz)\\
 &=& \int_{\Omega}\mid f(z(t_0)+s)\mid^2 \;q(dz),
\end{eqnarray*}
 it follows that $\{U_{s}\mid s\in\mathbb(R)\}$ is the unitary operator group on the Hilbert space $H$. In order to prove the following $L^1$ ergodic theorem, we need the following mean ergodic theorem:
\begin{theorem}   \textbf{(Mean ergodic theorem \cite{kre})}\label{t44}
 Let $\{U_{t}\mid t\in\mathbb(R)\}$ be the strong continuous one-parameter unitary operator group on the Hilbert space $H$. $P$ is the orthogonal project operator from $H$ to $\{f\in H\mid U_t f=f,t\in\mathbb{R}\}$, then
 \begin{equation}\label{e62}
s-\lim_{T\rightarrow \infty}\frac{1}{T}\int_0^T U_t \;dt=P.
\end{equation}

\end{theorem}
By which, we can prove the $L^1$ ergodic theorem of GIFs
\begin{theorem}   \textbf{($L^1$ ergodic theorem)}\label{t45}
 Let $(\Omega,q)$ be a GIF, then $\forall f(x)\in L^1(X)$, we have
\begin{equation}\label{e63}
\frac{1}{T}\int_0^T f(z(t+s))\; ds\stackrel{L^1_q}{\longrightarrow}f'(z),
\end{equation}
where $f'\in L^1_q(\Omega)$ and $\int_{\Omega}f'(z)\;q(dz)=\int_{\Omega}f'(\Gamma^s z)\;q(dz)$.

\end{theorem}
\noindent\textbf{Proof of Theorem \ref{t45}.} $\forall g(x)\in L^{\infty}_{\mu}(X)\subseteq L^{2}_{\mu}(X)$, by the mean ergodic theorem, we have
\begin{equation}\label{e64}
\frac{1}{T}\int_0^T g(z(t+s)) \;ds\stackrel{L^2_q}{\longrightarrow}g'(z)\in L^2_q(\Omega).
\end{equation}
Because
\begin{equation}\label{e65}
\parallel\frac{1}{T}\int_0^T g(z(t+s)) \;ds\parallel_{L^{\infty}_q(\Omega)}\leq \parallel g(x)\parallel_{L^{\infty}_{\mu}(X)},
\end{equation}
 we have
\begin{equation}\label{e66}
\mid\langle A^g_T,\chi_B\rangle\mid\leq\parallel g(x)\parallel_{L^{\infty}(X)}q(B), \forall B\in\mathcal{B}_{\Omega},
\end{equation}
here $A^g_T=\frac{1}{T}\int_0^T g(z(t+s)) ds$. Hence $\parallel g'\parallel_{L^{\infty}_q(\Omega)}\leq\parallel g\parallel_{L^{\infty}(X)}$, i.e. $g'\in L^{\infty}_q(\Omega)$. By $\parallel\cdot\parallel_{L^1_q}\leq\parallel\cdot\parallel_{L^2_q}$, we have
\begin{equation}\label{e67}
\frac{1}{T}\int_0^T g(z(t+s)) \;ds\stackrel{L^1_q}{\longrightarrow}g'(z)\in L^{\infty}_q(\Omega),
\end{equation}
so the theorem holds for the dense set $L^{\infty}_{\mu}(X)\subseteq L^1_{\mu}(X)$.
\par
$\forall f\in L^1_{\mu}(X)$, fix $\varepsilon>0$, we can choose $g\in L^{\infty}_{\mu}(X)$ with $\parallel f-g\parallel_{L^1_{\mu}(\Omega)}\leq\varepsilon$, then
\begin{equation}\label{e68}
\parallel\frac{1}{T}\int_0^T f(z(t+s))\; ds-\frac{1}{T}\int_0^T g(z(t+s)) \;ds\parallel_{L^1_q(\Omega)}\leq\varepsilon.
\end{equation}
We have proved that there exists $g'(z)\in L^{\infty}_q(\Omega)$ with
\begin{equation}\label{e69}
\parallel\frac{1}{T}\int_0^T g(z(t+s))\; ds-g'(z)\parallel_{L^1_q(\Omega)}\leq\varepsilon, for \ T>T_0.
\end{equation}
Hence we get
\begin{equation}\label{e70}
\parallel\frac{1}{T}\int_0^T f(z(t+s))\; ds-\frac{1}{T'}\int_0^{T'} f(z(t+s)) \; ds\parallel_{L^1_q(\Omega)}\leq4\varepsilon, for \ T,T'>T_0,
\end{equation}
which means that the ergodic averages generate a Cauchy sequence in $L^1_q(\Omega)$. So they have a limit $f'(z)\in L^1_q(\Omega)$. Because
\begin{equation}\label{e71}
\parallel\frac{1}{T}\int_0^T f(z(t+s))\; dt-\frac{1}{T}\int_0^{T} f(z(t))\; dt\parallel_{L^1_q(\Omega)}\leq\frac{2}{T}\parallel f(x)\parallel_{L^1_{\mu}(X)},
\end{equation}
we know that $f'(z)$ satisfies $\int_{\Omega}f'(z)\;q(dz)=\int_{\Omega}f'(\Gamma^s z)\;q(dz)$.
\begin{flushright}
  $\Box$
\end{flushright}
\par
Now we can prove the pointwise ergodic theorem of GIFs
\begin{theorem}   \textbf{(Pointwise ergodic theorem)}\label{t46}
 Let $(\Omega,q)$ be a GIF, $\phi(x)\in L^1(X)$. Then for $q-$almost all $z(t)\in\Omega$, we have
\begin{equation}\label{e72}
\lim_{T\rightarrow\infty}\frac{1}{T}\int_0^T \phi(z(t))\; dt=\Phi(z),
\end{equation}
here $\Phi\in L^1_q(\Omega)$ and $\int_{\Omega}\Phi(z)\;q(dz)=\int_{X}\phi(x)\;\mu(dx)$.

\end{theorem}
\noindent\textbf{Proof of Theorem \ref{t46}.} Assume that $\phi_0\in L^{\infty}_{\mu}(X)$. Then by theorem \ref{t45}, we know that
\begin{equation}\label{e73}
A_T(\phi_0)=\frac{1}{T}\int_{0}^{T}\phi_0(z(t))\;dt\stackrel{L^1_q}{\longrightarrow}\Phi_0\in L^1_q(\Omega), as \ T\rightarrow\infty,
\end{equation}
where $\Phi_0$ satisfies $\int_{\Omega}\Phi_0(z)\;q(dz)=\int_{\Omega}\Phi_0(\Gamma^s z)\;q(dz)$.\\
Then, for given $\varepsilon>0$, choose some $M$, such that
\begin{equation}\label{e74}
\parallel\Phi_0-A_M(\phi_0)\parallel_{L^1_q(\Omega)}<\varepsilon^2.
\end{equation}

We need a more general maximal ergodic theorem as what follows:
\begin{theorem}   \textbf{(Maximal ergodic theorem)}\label{t47}
 Let $(\Omega,q)$ be a GIF, for $f(z)\in L^1_q(\Omega)$ satisfies $\int_{\Omega}f(z)\;q(dz)=\int_{\Omega}f(\Gamma^s z)\;q(dz)$ and $\alpha >0$, we define the maximal function
 \begin{equation}\label{e411}
  f^{*}(z)=\sup_{T>0}\frac{1}{T}\int_{0}^{T}f(\Gamma^t z)\; dt,
\end{equation}
and the set
 \begin{equation}\label{e421}
  E_{\alpha}^{f}=\{z\in\Omega\mid f^{*}(z)>\alpha\},
\end{equation}
then we have
 \begin{equation}\label{e431}
  \alpha q(E^{f}_{\alpha})\leq\parallel f\parallel_{L^1_q(\Omega)}.
\end{equation}
\end{theorem}

The proof of this theorem is the same as that of theorem \ref{t41}.
Applying  theorem  \ref{t47} to the function $\psi(z)=\Phi_0-A_M(\phi_0)$,
 we get
\begin{equation}\label{e75}
\varepsilon q(\{z\in\Omega\mid\sup_{T>0}(A_T(\Phi_0-A_M(\phi_0)))\mid>\varepsilon\})<\varepsilon^2.
\end{equation}
Clearly, $A_T(\Phi_0)=\Phi_0$. And we have
\begin{eqnarray*}
A_T(A_M(\phi_0)) &=& \frac{1}{TM}\int_{0}^{T}\int_{0}^{M}\phi_0(z(t+s))\;dt\;ds\\
&=& A_T(\phi_0)+O_M(\frac{\parallel \phi_0\parallel_{\infty}}{T}).
\end{eqnarray*}
If $M$ is fixed and $T\rightarrow\infty$, we see that
\begin{eqnarray*}
q(\{z\in\Omega\mid\limsup_{T\rightarrow\infty}\mid\Phi_0-A_T(\phi_0)\mid>\varepsilon\})
&=& q(\{z\in\Omega\mid\limsup_{T\rightarrow\infty}\mid\Phi_0-A_T(A_M(\phi_0))\mid>\varepsilon\})\\
&\leq& q(\{z\in\Omega\mid\sup_{T>0}(A_T(\Phi_0-A_M(\phi_0)))\mid>\varepsilon\})<\varepsilon,
\end{eqnarray*}
 which means that $A_T(\phi_0)\rightarrow \Phi_0$ for $q-$almost all $z$.
\par
For any $\phi\in L^1_{\mu}(X)$, fix $\varepsilon>0$ and choose $\phi_0\in L^{\infty}_{\mu}(X)$ with $\parallel \phi-\phi_0\parallel_{L^1_{\mu}(X)}\leq\varepsilon^2$. Let $\Phi\in L^1_q(\Omega)$ be the $L^1_{q}-$limit of $A_T(\phi)=\frac{1}{T}\int_{0}^{T}\phi(z(t))\;dt$ and $\Phi_0\in L^1_q(\Omega)$ be the $L^1_{q}-$limit of $A_T(\phi_0)=\frac{1}{T}\int_{0}^{T}\phi_0(z(t))\;dt$. Because $\parallel A_T(\phi)-A_T(\phi_0)\parallel_{L^1_q(\Omega)}\leq\parallel \phi-\phi_0\parallel_{L^1_{\mu}(X)}$, we have $\parallel \Phi-\Phi_0\parallel_{L^1_q(\Omega)}\leq\varepsilon^2$. From this and the maximal ergodic theorem \ref{t41}, we get
\begin{equation*}
\begin{split}
&\quad\; q(\{z\in\Omega\mid\limsup_{T\rightarrow\infty}\mid\Phi-A_T(\phi)\mid>2\varepsilon\})\\
&\leq q(\{z\in\Omega\mid(\mid\Phi-\Phi_0\mid+\limsup_{T\rightarrow\infty}\mid\Phi_0-A_T(\phi_0)\mid+\sup_{T>0}\mid A_T(\phi-\phi_0)\mid)>2\varepsilon\})\\
&\leq q(\{z\in\Omega\mid\Phi-\Phi_0\mid>\epsilon\})+q(\{z\in\Omega\mid\sup_{T>0}(A_T\mid(\phi-\phi_0)\mid>\varepsilon\})\\
&\leq\varepsilon^{-1}\parallel\Phi-\Phi_0\parallel_{L^1_q(\Omega)}+\varepsilon^{-1}\parallel\phi-\phi_0\parallel_{L^1_{\mu}(X)}\leq2\varepsilon,
\end{split}
\end{equation*}
which shows that $A_T(\phi)\rightarrow \Phi$ for $q-$almost all $z$. And it is clearly that $\int_{\Omega}\Phi(z)\;q(dz)=\int_{X}\phi(x)\;\mu(dz)$ by the Fubini's theorem.
\begin{flushright}
  $\Box$
\end{flushright}

\section{More on Generalized incompressible flows }
\par
Let $X$ be a compact metric space. We know that $\Omega=X^{\mathbb{R}}$ is compact with respect to the product topology, i.e. the pointwise convergence topology, but we can not find a metric compatible to this topology. On the other hand, we may define a metric on $\Omega$ as follows
\begin{equation}\label{e76}
  d(z(t),z'(t))=\left\{
     \begin{array}{ll}
       \sup_{t\in\mathbb{R}}\mid z(t)-z'(t)\mid & ,if \mid z(t)-z'(t)\mid<1,\forall t \\
       1 & ,otherwise.
     \end{array}
   \right.
\end{equation}
Hence, we get a uniformly convergence topology, and $\Omega$ is a metric space. Although we cannot make $\Omega$ compact and metrizable at the same time, we still can get a compact Hausdorff space because the product space of Hausdorff spaces is still Hausdorff with respect to the product topology.
\par
Then using the Riesz representation theorem on compact Hausdorff space and the Krylov-Bogolioubov theorem, we can prove that
\begin{proposition}   \label{p51}
 The set of all GIFs denoted by $\text{GIF}\,(\Omega)$ is non-empty.
\end{proposition}
\noindent\textbf{Proof.} Recall that a generalized flow is a probability measure on the path space $\Omega=X^\mathbb{R}$. Indeed, map $\Gamma^s:\Omega\rightarrow\Omega$ introduces a map on the set of generalized flows $\text{GF}\,(\Omega)$:
\begin{eqnarray*}\label{e77}
  (\Gamma_{s})^{\sharp}: &&\text{GF}\,(\Omega)\longrightarrow \text{GF}\,(\Omega) \\
   &&q\longmapsto (\Gamma_{s})^{\sharp}q,
\end{eqnarray*}
where $(\Gamma_{s})^{\sharp}q(B)=q((\Gamma^s)^{-1}B),\; \forall B\in\mathcal{B}_{\Omega}$. If $\{\omega_n\}_{n=1}^{\infty}$ is a sequence of generalized flows, we can construct a new sequence $\{q_n\}_{n=1}^{\infty}$ by
\begin{equation}\label{e78}
  q_n=\frac{1}{n}\int_0^n(\Gamma_{s})^{\sharp}\omega_n\; ds.
\end{equation}
By the weak$-^{*}$ compactness of $\text{GF}\,(\Omega)$, there exists $q$ such that a subsequence  $q_{n_j}$ weak$-^{*}$ converge to $q$. We can prove that $q\in \text{GIF}\,(\Omega)$, even more, $q$ is an invariant measure of $\Gamma^s$ on $\Omega$, which is a smaller subset of $\text{GIF}\,(\Omega)$ as we mentioned before.
\begin{equation*}\label{e79}
  \begin{split}
  &\quad\; \mid\int_{\Omega}f(\Gamma^s z)\;q(dz)-\int_{\Omega}f(z)\;q(dz)\mid\\
  &=\lim_{j\rightarrow\infty}\mid\int_{\Omega}f(\Gamma^s z)\;q_{n_j}(dz)-\int_{\Omega}f(z)\;q_{n_j}(dz)\mid\\
  &=\lim_{j\rightarrow\infty}\mid\frac{1}{n_j}\int_0^{n_j}\int_{\Omega}f(\Gamma^s z)(\Gamma^t)^{\sharp}\omega_{n_j}(dz)\;dt
  -\frac{1}{n_j}\int_0^{n_j}\int_{\Omega}f(z)(\Gamma^t)^{\sharp}\omega_{n_j}(dz)\;dt\mid\\
  &=\lim_{j\rightarrow\infty}\frac{1}{n_j}\mid\int_s^{n_j+s}\int_{\Omega}f(\Gamma^t z)\omega_{n_j}(dz)\;dt-\int_0^{n_j}\int_{\Omega}f(\Gamma^t z)\omega_{n_j}(dz)\;dt\mid\\
  &=\lim_{j\rightarrow\infty}\frac{1}{n_j}\mid\int_{n_j}^{n_j+s}\int_{\Omega}f(\Gamma^t z)\omega_{n_j}(dz)\;dt-\int_0^{s}\int_{\Omega}f(\Gamma^t z)\omega_{n_j}(dz)\;dt\mid\\
  &\leq\lim_{j\rightarrow\infty}\frac{1}{n_j}2s\parallel f\parallel_{\infty}=0,
  \end{split}
\end{equation*}
so we have $\int_{\Omega}f(\Gamma^s z)\;q(dz)=\int_{\Omega}f(z)\;q(dz)$, i.e. $q\in \text{GIF}\,(\Omega)$.
\begin{flushright}
  $\Box$
\end{flushright}
\par
Before continuing our discussion, we define the weak ergodicity
\begin{definition}   \textbf{(Weak Ergodic GIF)}\label{d51}
Let $(\Omega,q)$ be a GIF and set $\Phi=\{z\in\Omega\mid z(0)\in E,\;E\in\mathcal{B}_X\}$ be the set of all paths which is in $E$ when $t=0$. We call $q$ is weak ergodic, if for any set $\Phi$ of above type and also satisfies $\Gamma^s\Phi=\Phi$, $\forall s\in\mathbb{R}$, we have $q(\Phi)=0$ or 1.
\end{definition}
In section 3, we have proved that the GIF $q$ induced from a classical ergodic flow is ergodic. Actually, we can prove that all the ergodic GIF can be defined in a similar way.
\begin{proposition}   \label{p52}
 If the GIF $q$ is ergodic, then $q-$almost all $z\in\Omega$ is ergodic.
\end{proposition}
\noindent\textbf{Proof.} Choose $\hat{z}\in\Omega$, we can see that set $\Psi=\{\Gamma^s \hat{z}\mid s\in\mathbb{R}\}$ is invariant under $\Gamma^s$. Hence by the ergodicity of $q$, we have $q(\Psi)=0$ or 1, here $q(\Psi)=1$ means that $q$ is supported on $\Psi$.
For any fixed path $z\in\Omega=X^\mathbb{R}$, there is a natural time-dependent map on the subset of $X$ induced by $z$, we define a flow $T^t$ on $\{z(t)\mid t\in\mathbb{R}\}\subseteq X$ by
 $$T^t:z(t_0)\mapsto z(t_0+t).$$
 Let $B\subseteq X$ be an invariant set  of flow $\hat{T}^t$ induced by $\hat{z}$. Define $\Phi=\{\hat{z}\in\Omega\mid \hat{z}(0)\in B\}$. Because $B$ is invariant under $\hat{T}^t$, we get $\Gamma^{-s}\Phi=\Phi,\forall s\in\mathbb{R}$, so $q(\Phi)=0$ or 1. Then by the incompressible condition, we have
 \begin{eqnarray*}\label{e80}
  \mu(B)&=&\int_X\chi_B\;\mu(dx)=\int_{\Omega}\chi_{B}(z(0))\;q(dz)\\
  &=&\int_{\Omega}\chi_{\Phi}(z)\;q(dz)=q(\Phi)=0 \;or\; 1,
\end{eqnarray*}
which shows that the flow $\hat{z}(t)$ is ergodic, so is $\Gamma^s\;\hat{z}(t)$, for any $s\in\mathbb{R}$, which means that $q-$almost all $z\in\Omega$ is ergodic.
\begin{flushright}
  $\Box$
\end{flushright}
\begin{remark}   \label{r51}
 (a) Choose $z(t)\in supp\;q$, by Proposition \ref{p52}, $z(t)$ is ergodic. So $T^t$ is defined on $X$ except a zero measure set. We have
\begin{equation}\label{e81}
  \int_{\Omega}f(z(t))\;q(dz)=\int_{X}f(x\mapsto T^t(x))\;\mu(dx), \forall f(x)\in C(X),
\end{equation}
which is similar to (\ref{e31}).
\\
(b) Weak ergodic GIFs is a larger class than ergodic GIFs. We can see that for an ergodic GIF, it can only support on $\{\Gamma^s\;z\mid s\in\mathbb{R}\}$, here $z$ is a classical ergodic flow. However, for weak ergodic GIF, it can support on $\{\Gamma^s\;x\mid s\in\mathbb{R}\}\cup\{\Gamma^s\;y\mid s\in\mathbb{R}\}$, if $x(t)$ and $y(t)$ are ergodic flows with respect to the same measure $\mu$ on $X$.
 \end{remark}
Before giving a theorem about some properties of $\text{GIF}\,(\Omega)$, we prove a useful lemma
\begin{lemma} \label{t53}
Let $(\Omega,q)$ be a GIF, and for any $E\in\mathcal{B}_X$, define $\Phi=\{z\in\Omega\mid z(0)\in E\}$, then $q(\Gamma^{-s}\Phi)=q(\Phi)$.
\end{lemma}
\noindent\textbf{Proof.} we have
 \begin{eqnarray*}\label{}
 q(\Gamma^{-s}\Phi)&=&\int_{\Omega}\chi_{\Gamma^{-s}\Phi}(z)\;q(dz)=\int_{\Omega}\chi_{E}(z(s))\;q(dz)\\
 &=&\int_{\Omega}\chi_{E}(z(0))\;q(dz)=\int_{\Omega}\chi_{\Phi}(z)\;q(dz)=q(\Phi),\\
\end{eqnarray*}
which proves the lemma.
\begin{flushright}
  $\Box$
\end{flushright}

Now we prove a theorem about the properties of the set $\text{GIF}\,(\Omega)$.
\begin{theorem}   \label{t54}
 (1) $\text{GIF}\,(\Omega)$ is a compact subset of $\text{GF}\,(\Omega)$.\\
 (2) $\text{GIF}\,(\Omega)$ is convex.\\
 (3) $q$ is an extreme point of $\text{GIF}\,(\Omega)$ iff $q$ is a weak ergodic GIF.
\end{theorem}
\noindent\textbf{Proof.}(1) Suppose $\{q_n\}_1^{\infty}$ is a sequence of members of $\text{GIF}\,(\Omega)$ and $q_n\rightarrow q$ in $\text{GF}\,(\Omega)$. Then
 \begin{eqnarray*}\label{e82}
  \int_{\Omega}f(z(t_0))\;q(dz)&=&\lim_{n\rightarrow{\infty}}\int_{\Omega}f(z(t_0))\;q_n(dz)\\
  &=&\lim_{n\rightarrow{\infty}}\int_{\Omega}f(z(0))\;q_n(dz)\\
  &=&\int_{\Omega}f(z(0))\;q(dz),
\end{eqnarray*}
we have $q\in \text{GIF}\,(\Omega)$. Hence $\text{GIF}\,(\Omega)$ is a closed subset of compact set $\text{GIF}\,(\Omega)$. So it must be compact.\\
(2) It is clear that if $q_1,q_2\in \text{GIF}\,(\Omega)$, then $(1-s)q_1+sq_2$ is a GIF. \\
(3) $(\Rightarrow)$ Suppose $q\in \text{GIF}\,(\Omega)$ and it is not weak ergodic, then there exists a Borel set $E\subset X$. We have $\Phi=\{z\in\Omega\mid z(0)\in E\}$ such that $\Gamma^{-s}\Phi=\Phi, \forall s\in\mathbb{R}$ with $0<q(\Phi)<1$.\\
Define $q_1$ and $q_2$ by
 \begin{eqnarray}\label{e83}
  q_1(B)&=&\frac{q(B\cap\Phi)}{q(\Phi)},\\
  q_2(B)&=&\frac{q(B\cap(\Omega\setminus\Phi)}{q(\Phi)}.
\end{eqnarray}
Next we prove that $q_1,\;q_2\in\text{GIF}\,(\Omega)$. To prove $q_1,\;q_2$ are incompressible, it suffices to prove that, for any $Q\subset X$,
  \begin{eqnarray}\label{e84}
  \int_{\Omega}\chi_Q(z(t_0))q_1(dz)&=&\int_{\Omega}\chi_Q(z(0))q_1(dz),\\
  \int_{\Omega}\chi_Q(z(t_0))q_2(dz)&=&\int_{\Omega}\chi_Q(z(0))q_2(dz).
\end{eqnarray}
Also, notice that
\begin{equation}\label{e85}
  q_1(\Psi(t_0))=\int_{\Omega}\chi_Q(z(t_0))q_1(dz)=\int_{\Omega}\chi_{\Psi(t_0)}(z)q_1(dz),
\end{equation}
where $\Psi(t_0)=\{z\in\Omega\mid z(t_0)\in Q\}$.\\
Because $\Gamma^{-s}\Phi=\Phi$, we get
\begin{equation}\label{e86}
 \Phi=\Phi(t_0)=\{z\in\Omega\mid z(t_0)\in E\},\forall t_0\in\mathbb{R}.
\end{equation}
Hence,
\begin{equation}\label{e87}
 q_1(\Psi(t_0))=\frac{q(\Psi(t_0)\cap\Phi)}{q(\Phi)}=\frac{q(\Psi(t_0)\cap\Phi(t_0))}{q(\Phi)}.
\end{equation}
Because $\Psi(t_0)\cap\Phi(t_0)=\{z\in\Omega\mid z(t_0)\in Q\cap E\}=\Gamma^{-t_0}(\Psi(0)\cap\Phi(0))$, and by lemma \ref{t53}, we have
\begin{equation}\label{e88}
\frac{q(\Psi(t_0)\cap\Phi(t_0))}{q(\Phi)}=\frac{q(\Psi(0)\cap\Phi(0))}{q(\Phi)}=q(\Psi(0)).
\end{equation}
Combining (\ref{e85}), (\ref{e87}), (\ref{e88}), we get (\ref{e84}) which proves $q_1\in \text{GIF}\,(\Omega)$. Similarly, we can prove $q_2\in\text{GIF}\,(\Omega)$. So we have $q_1\neq q_2\in\text{GIF}\,(\Omega)$ such that
\begin{equation}\label{e89}
q(B)=q(\Phi)q_1(B)+(1-q(\Phi))q_2(B).
\end{equation}
$(\Leftarrow)$ Suppose  $\text{GIF}\,q$ is weak ergodic, and
\begin{equation}\label{e90}
q(B)=pq_1(B)+(1-p)q_2(B),
\end{equation}
where $q_1,q_2\in \text{GIF}\,(\Omega)$ and $0<p<1$. We will prove $q_1=q_2$. As we showed in remark \ref{r21}, we can find $\mu_1,\mu_2$ such that
\begin{equation}\label{e91}
  \int_{\Omega}f(z(t_{0}))\;q_i(dz)=\int_{\Omega}f(z(0))\;q_i(dz)=\int_{X}f(x)\;\mu_i(dx),\; \forall f(x)\in C(X), \; \forall t_{0}\in \mathbb{R},\; i=1,2,
\end{equation}
Clearly, $q_1\ll q$ and $\mu_1\ll\mu$. Let $E\in\mathcal{B}_X$, $\Phi=\{z\in\Omega\mid z(0)\in E\}$. We have
 \begin{eqnarray*}\label{e92}
  q_1(\Phi)&=&\mu_1(E)=\int_{E}\frac{d\mu_1}{d\mu}(x)\;\mu(dx)\\
  &=&\int_{\Omega}\chi_{E}(z(0))\frac{d\mu_1}{d\mu}(z(0))\;q(dz)\\
  &=&\int_{\Omega}\chi_{\Phi}(z)\frac{d\mu_1}{d\mu}(z(0))\;q(dz)\\
  &=&\int_{\Phi}\frac{d\mu_1}{d\mu}(z(0))\;q(dz)=\int_{\Phi}\frac{dq_1}{dq}(z)\;q(dz).\\
\end{eqnarray*}
By Radon-Nikodym theorem, for $q-$almost everywhere on $\Phi$, we have $\frac{d\mu_1}{d\mu}(z(0))=\frac{dq_1}{dq}(z)$.
\par
Let $E=\{x\mid \frac{d\mu_1}{d\mu}(x)<1\}$, then $\Phi=\{z\in\Omega\mid z(0)\in E\}=\{z\in\Omega\mid \frac{dq_1}{dq}(z)<1\}$.
We have
  \begin{equation*}\label{e93}
 \begin{split}
 &\quad\;\int_{\Phi\cap\Gamma^{-s}\Phi} \frac{dq_1}{dq}(z)\;q(dz)+\int_{\Phi\setminus\Gamma^{-s}\Phi}\frac{dq_1}{dq}(z)\;q(dz)\\
 &=q_1(\Phi)=q_1(\Gamma^{-s}\Phi)\\
 &=\int_{\Phi\cap\Gamma^{-s}\Phi} \frac{dq_1}{dq}(z)\;q(dz)+\int_{\Gamma^{-s}\Phi\setminus\Phi}\frac{dq_1}{dq}(z)\;q(dz),
 \end{split}
\end{equation*}
here we use lemma \ref{t53}.
\par
Since $\frac{dq_1}{dq}<1$ on $\Phi\setminus\Gamma^{-s}\Phi$ and $\frac{dq_1}{dq}\geq1$ on $\Gamma^{-s}\Phi\setminus\Phi$, also we have
\begin{eqnarray*}\label{e94}
 q(\Gamma^{-s}\Phi\setminus\Phi)&=&q(\Gamma^{-s}\Phi)-q(\Phi\cap\Gamma^{-s}\Phi)\\
 &=&q(\Phi)-q(\Phi\cap\Gamma^{-s}\Phi)\\
 &=&q(\Phi\setminus\Gamma^{-s}\Phi),
\end{eqnarray*}
so $q(\Gamma^{-s}\Phi\setminus\Phi)=q(\Phi\setminus\Gamma^{-s}\Phi)=0$, i.e. $q(\Gamma^{-s}\Phi\vartriangle\Phi)=0$. Therefore, by the weak ergodicity of $q$, we know $q(\Phi)=0$ or 1. But if $q(\Phi)=1$, then $q(\Omega)=\int_{\Phi}\frac{dq_1}{dq}(z)q(dz)<q(\Phi)=1$, this is contradict to $q(\Omega)=1$, hence $q(\Phi)=0$.
\par
Similarly, we can prove that $q(\{z\in\Omega\mid\frac{dq_1}{dq}(z)>1\})=0$. So for $q-$almost everywhere we have  $q_1=q$. Therefore we prove that $q$ is an extreme point of $\text{GIF}\,(\Omega)$.
\begin{flushright}
  $\Box$
\end{flushright}
\begin{remark}\label{r52}
From this theorem, we can see that  the extreme points of $\text{GIF}\,(\Omega)$ are not ergodic but weak ergodic,  unlike the classical result of the set of invariant measures for continuous transformations.
\end{remark}




  \end{document}